\documentclass[11pt]{article}
\usepackage{CJK}
\usepackage{amsmath}
\usepackage{amsmath,amssymb,latexsym,color}
\usepackage[mathscr]{eucal}
\usepackage{graphicx}
\usepackage{pstricks}
\newtheorem{thm}{Theorem}[section]
\newtheorem{defin}[thm]{Definition}
\newtheorem{prop}[thm]{Proposition}

\newtheorem{cor}[thm]{Corollary}

\newtheorem{example}{Example}

\newcommand{\proof}{{\it Proof.\quad}}
\newcommand{\qed}{\hfill\Box\medskip}
\begin{document}
\begin{CJK*}{GBK}{song}

\title{\bf On the metric dimension of line graphs}

\author{Min Feng\quad Min Xu\quad Kaishun Wang\footnote{Corresponding author. E-mail address: wangks@bnu.edu.cn}\\
{\footnotesize   \em  Sch. Math. Sci. {\rm \&} Lab. Math. Com.
Sys., Beijing Normal University, Beijing, 100875,  China}}
\date{}
\maketitle

\begin{abstract}
Let $G$ be a (di)graph. A set $W$ of vertices in $G$ is a
\emph{resolving set} of $G$ if every vertex $u$ of $G$ is uniquely
determined by its vector of distances to all the vertices in $W$.
The \emph{metric dimension} $\mu (G)$ of $G$ is the minimum
cardinality of all the resolving sets of $G$.  C\'aceres et al.
\cite{Ca2} computed the metric dimension of the line graphs of
complete bipartite graphs. Recently, Bailey and Cameron \cite{Ba}
computed the metric dimension of the line graphs of complete
graphs. In this paper we study the metric dimension of the line
graph $L(G)$ of $G$. In particular, we show that
$\mu(L(G))=|E(G)|-|V(G)|$ for a strongly connected digraph $G$
except for directed cycles, where $V(G)$ is the vertex set and
$E(G)$ is the edge set of $G$. As a corollary, the metric
dimension of de Brujin digraphs and Kautz digraphs is given.
Moreover, we prove that
$\lceil\log_2\Delta(G)\rceil\leq\mu(L(G))\leq |V(G)|-2$ for a
simple connected graph $G$ with at least five vertices, where
$\Delta(G)$ is the maximum degree of $G$. Finally, we obtain the
metric dimension of the line graph of a tree in terms of its
parameters.


%\medskip
%\noindent {\em 2010 MSC:} 05C12, 05E30.
\end{abstract}

\section{Introduction}

Let $G$ be a (di)graph. We often write $V(G)$ for the vertex set
of $G$ and $E(G)$ for the edge set of $G$. A (di)graph $G$ is
(strongly) connected if for any two distinct vertices $u$ and $v$
of $G$, there exists a path from $u$ to $v$. In this paper we only
consider finite strongly  connected digraphs, or undirected simple
connected graphs. For two vertices $u$ and $v$ of $G$, we denote
the distance from $u$ to $v$ by $d_{G}(u,v)$. A \emph{resolving
set} of $G$ is a set of vertices $W=\{w_{1},\ldots,w_{m}\}$ such
that for each $u\in V(G)$, the vector $D(u|W) =
(d_G(u,w_{1}),\ldots,d_G(u,w_{m}))$
    uniquely determines $u$. The
\emph{metric dimension} of $G$, denoted by $\mu(G)$, is the minimum
cardinality of all the resolving sets of $G$.

Metric dimension of graphs was introduced in the 1970s,
independently by Harary and Melter \cite{Ha} and by Slater
\cite{Sl1}. Metric dimension of digraphs was first studied by
Chartrand et al. in \cite{Char1} and further in \cite{Char2}. Fehr
et al. \cite{Fe} investigated  the metric dimension of Cayley
digraphs. In graph theory, metric dimension is a parameter that
has appeared in various applications, as diverse as   network
discovery and verification \cite{Bee}, strategies for the
Mastermind game \cite{Chv},  combinatorial optimization \cite{Se}
and so on. It was noted in \cite[p. 204]{Ga} and \cite{Kh} that
 determining the metric dimension of a graph is an NP-complete problem.

Let $L(G)$ denote the line graph of a (di)graph $G$. For the
complete bipartite graph $K_{m,n}$, C\'aceres et al. \cite{Ca2}
proved that
\begin{equation*}
\mu(L(K_{m,n}))=\left\{
\begin{array}{ll}
\lfloor\frac{2(m+n-1)}{3}\rfloor, &m\leq n\leq 2m-1,n\geq 2,\\
n-1,&n\geq 2m.
\end{array}\right.
\end{equation*}
  For the
complete graph $K_n$ when $n\geq 6$, Bailey and Cameron \cite{Ba}
proved that $\mu(L(K_n))=\lceil\frac{2n}{3}\rceil$.

Motivated by these results, in this paper  we study the metric
dimension of the line graph of a (di)graph. In Section 2, we show
that $\mu(L(G))=|E(G)|-|V(G)|$ for a strongly connected  digraph
$G$ except for directed cycles. As a corollary, the metric
dimension of de Brujin digraphs and Kautz digraphs, which are two
families of famous networks, is given. In Section 3, we prove that
$\lceil\log_2\Delta(G)\rceil\leq\mu(L(G))\leq |V(G)|-2$ for a
connected graph $G$ with at least five vertices, where $\Delta(G)$
is the maximum degree of $G$. Finally, we obtain the metric
dimension of the line graph of a tree in terms of its parameters.

\section{Line graph of a digraph}

Let $G$ be a digraph.  For a directed edge $a=(x,y)$ of $G$, we
say that $x$ is the \emph{head} of $a$ and $y$ is the \emph{tail}
of $a$; we also say that $a$ is the \emph{out-going edge} of $x$
and the \emph{in-coming edge} of $y$. For $x\in V(G)$, we denote
the set of all out-going edges of $x$ by $E_{G}^{+}(x)$ and the
set of all in-coming edges of $x$ by $E_{G}^{-}(x)$. The
\emph{line  graph} of $G$ is the digraph $L(G)$ with the edges of
$G$ as its vertices, and where $(a,b)$ is a directed edge in
$L(G)$ if and only if the tail of $a$ is the head of $b$ in $G$.
For two distinct vertices $a=(x_1,x_2), b=(y_1,y_2)$ of $L(G)$, we
have
\begin{equation}\label{distance}
d_{L(G)}(a,b)=d_{G}(x_2,y_1)+1.
\end{equation}

Note that $\mu(L(G))=1$ if $G$ is a directed cycle.

\begin{thm}\label{main}
If $G$ is a strongly connected digraph except for  directed
cycles, then
$$\mu(L(G))=|E(G)|-|V(G)|.$$
\end{thm}
\proof
 Let $R$ be a resolving set of $L(G)$ with the minimum
cardinality. For each vertex $x$ of $G$, since $G$ is strongly
connected, $E^-_G(x)\neq\emptyset$. If $|E_G^-(x)|\geq2$, pick two
distinct edges $a,b\in E^-_G(x)$. For any $c\in
V(L(G))\backslash\{a,b\}$, since $d_{L(G)}(a,c)=d_{L(G)}(b,c)$,
$a\in R$ or $b\in R$. It follows that $|E_G^-(x)\cap
R|\geq|E_G^-(x)|-1$. If $|E_G^-(x)|=1$, the above inequality is
directed. By $R=\dot{\cup}_{x\in V(G)}(E_G^-(x)\cap R),$ we obtain
\begin{equation}\label{low}
\mu(L(G))=|R|\geq\sum_{x\in V(G)}(|E^-_{G}(x)|-1)=|E(G)|-|V(G)|.
\end{equation}

 Let $W$
be a set obtained from $E(G)$ by deleting one in-coming edge of
each vertex of $G$. Since $G$ is not a directed cycle,
$W\neq\emptyset$. We shall prove that $W$ is a resolving set of
$L(G)$. It suffices to show that, for any two distinct edges
$a=(x_1,x_2)$ and $b=(y_1,y_2)$ in $E(G)\backslash W$,
 there exists an edge $c\in W$ such that
\begin{equation}\label{ineq}
d_{L(G)}(a,c)\neq d_{L(G)}(b,c).
\end{equation}
Let $A$ denote the set of all the heads of each edge of $W$. Pick
$z_0\in A$ satisfying $d_{G}(x_2,z_0)\leq d_G(x_2,z)$ for any
$z\in A$.

\em Case 1\rm. $d_G(x_2,z_0)\neq d_G(y_2,z_0)$. Pick $c\in
E^+_G(z_0)\cap W$. By (\ref{distance}), (\ref{ineq}) holds.

\em Case 2\rm. $d_G(x_2,z_0)=d_G(y_2,z_0)$. Owing to $a,b\not\in
W$, $x_2\neq y_2$, which implies $z_0\neq x_2$. Let
$P_{x_2,z_0}=(v_0=x_2,v_1,\ldots,v_k=z_0)$ be a shortest path from
$x_2$ to $z_0$ and $P_{y_2,z_0}=(u_0=y_2,u_1,\ldots,u_k=z_0)$ be a
shortest path from $y_2$ to $z_0$. Suppose $i$ denotes the minimum
index such that $v_i=u_i$. Since
$d_G(x_2,v_{i-1})<d_G(x_2,v_i)\leq d_G(x_2,z_0)$, we have
$v_{i-1}\not\in A$, which implies  $(v_{i-1},v_{i})\not\in W$.
Hence $(u_{i-1},u_i)\in W$ and $u_{i-1}\in A$. Pick
$c=(u_{i-1},u_{i})$. By (\ref{distance}), we have
\begin{equation*}
\begin{array}{ccl}
d_{L(G)}(a,c)&=   &d_G(x_2,u_{i-1})+1\\
             &\geq&d_G(x_2,z_0)+1\\
             &=   &d_G(y_2,z_0)+1\\
             &\geq&d_G(y_2,u_i)+1\\
             &=   &d_{L(G)}(b,c)+1\\
             &>   &d_{L(G)}(b,c),
\end{array}
\end{equation*}
so (\ref{ineq}) holds.

Therefore, $W$ is a resolving set of $L(G)$ with size
$|E(G)|-|V(G)|$, which implies that $\mu(L(G))\leq|E(G)|-|V(G)|$.
By (\ref{low}), the desired result follows.$\qed$

Let $K_d$ be the complete digraph with $d$ vertices. A
\emph{flowered complete digraph} of order $d$, denoted by $K_d^+$,
is a digraph obtained from $K_d$ by appending a self-loop at each
vertex. Let
\begin{equation*}
\begin{array}{c}
B(d,1)=K^+_d,\;B(d,n)=L(B(d,n-1));\\
K(d,1)=K_{d+1},\;K(d,n)=L(K(d,n-1)).
\end{array}
\end{equation*}
Then $B(d,n)$ is the \emph{de Brujin digraph} and $K(d,n)$ is the
\emph{Kautz digraph}. By \cite[Chapter 3]{Xu}, $B(d,n)$ and
$K(d,n)$ are strongly connected and
\begin{equation*}
\begin{array}{c}
|V(B(d,n))|=d^{n},\;|E(B(d,n))|=d^{n+1};\\
|V(K(d,n))|=d^{n}+d^{n-1},\;|E(K(d,n))|=d^{n+1}+d^{n}.
\end{array}
\end{equation*}
As a corollary of Theorem \ref{main}, we get the metric dimension of
de Brujin digraphs and Kautz digraphs, respectively.

\begin{cor}\label{dBK}
Let integers $d\geq 2$ and $n\geq 1$. Then

{\rm(i)} $\mu(B(d,n))=d^{n-1}(d-1)$;

{\rm(ii)} $\mu(K(d,n))=\left\{
\begin{array}{ll}
d,&\textup{ if }n=1,\\
d^{n-2}(d^{2}-1),&\textup{ if }n\geq 2.
\end{array}\right.$
\end{cor}

\section{Line graph of a graph}

Let $G$ be a graph with at least two vertices. The \emph{line
graph} of $G$ is the graph $L(G)$ with the edges of $G$ as its
vertices, and where two edges of $G$ are adjacent in $L(G)$ if and
only if they are adjacent in $G$.

If $G$ has at most four vertices, it is routine to compute the
metric dimension of $L(G)$. Next we shall consider the case
$|V(G)|\geq 5$.

\begin{thm}\label{upper}
If $G$ is  a connected graph with at least five vertices, then
$$\lceil\log_2\Delta(G)\rceil\leq\mu(L(G))\leq |V(G)|-2,$$ where
$\Delta(G)$ is the maximum degree of $G$.
\end{thm}
\proof Let $v$ be a vertex of degree $\Delta(G)$, and let
$\{f_1,\ldots,f_{\Delta(G)}\}$ be the set of all the edges
incident to $v$. Suppose $W=\{e_1,\ldots,e_{\mu(L(G))}\}$ is a
resolving set of $L(G)$ with the minimum cardinality. For each
$j\in\{1,\ldots,\mu(L(G)\}$, let $d_j=\min\{d_G(v,w)|w\textup{ is
incident to } e_j\}$. Then $d_{L(G)}(f_i,e_j)$ is  $d_j$ or
$d_j+1$. Therefore, the size of
$\mathcal{D}=\{D(f_i|W)|i=1,\ldots,\Delta(G)\}$ is at most
$2^{\mu(L(G))}$. Since $D(f_i|W)\neq D(f_k|W)$ for $i\neq k$,
$\Delta(G)\leq2^{\mu(L(G))}$, which implies the lower bound.

 Suppose $|V(G)|=5$. If $G$ is isomorphic to the path $P_5$
or the cycle $C_5$, since $\mu(L(P_5))=1$ and $\mu(L(C_5))=2$, the
upper bound is directed. If $G$ is not isomorphic to $P_5$ or
$C_5$, then $G$ has a subgraph $S$ isomorphic to $K_{1,3}$. Since
$E(S)$ is a resolving set of $L(G)$, $\mu(L(G))\leq 3$, which
implies the upper bound.

Now suppose $|V(G)|\geq 6$. Let $T$ be a spanning tree of $G$, and
let $v$ be a vertex of degree 1 in $T$. Suppose $T_1$ is the
subgraph of $T$ induced on $V(T)\backslash\{v\}$. We shall prove
that $E(T_1)$ is a resolving set of $L(G)$. It suffices to show
that, for any two distinct edges $a,b\in E(G)\backslash E(T_1)$,
there exists an edge $e\in E(T_1)$  such that
\begin{equation}\label{ineq1}
d_{L(G)}(a,e)\neq d_{L(G)}(b,e).
\end{equation}

{\em Case 1}. $a$ or $b$ is not incident to $v$. Without loss of
generality, suppose $a$ is not incident to $v$. Let $a=uu'$. Then
there exists a unique path $P_{u,u'}=(u_0=u,u_1,\ldots,u_k=u')$
between $u$ and $u'$ in $T$ where $k\geq 2$. If $b$ is not
adjacent to $u_0u_1$, then (\ref{ineq1}) holds for $e=u_0u_1\in
E(T_1)$; If $b$ is not adjacent to $u_{k-1}u_k$, then
(\ref{ineq1}) holds for $e=u_{k-1}u_k\in E(T_1)$. Now we assume
that $b$ is adjacent to both $u_0u_1$ and $u_{k-1}u_k$.

{\em Case 1.1}. $k=2$. Then $b$ is incident to $u_1$. Suppose
$b=u_1x$, where $x\in V(G)\backslash \{u_0,u_1,u_2\}$. Let
$S=\{u_0,u_1,u_2,x\}$ and $\overline{S}=V(T_1)\backslash S$. Since
$|V(T_1)|=|V(G)|-1\geq 5$, there exists an edge
$e\in[S,\overline{S}]_{T_1}$, where $[S,\overline{S}]_{T_1}$ is
the set of edges between $S$ and $\overline{S}$ in $T_1$. If $e$
is incident to $u_0$ or $u_2$, then $d_{L(G)}(a,e)=1$ and
$d_{L(G)}(b,e)=2$; If $e$ is incident to $u_1$ or $x$, then
$d_{L(G)}(a,e)=2$ and $d_{L(G)}(b,e)=1$. So (\ref{ineq1}) holds.

{\em Case 1.2}. $k\geq 3$. Note that $b$ is incident to $u_1$ or
$u_{k-1}$. Without loss of generality, assume that $b$ is incident
to $u_1$. Let $e=u_1u_2\in E(T_1)$. Then $d_{L(G)}(a,e)=2\neq
1=d_{L(G)}(b,e)$, (\ref{ineq1}) holds.

{\em Case 2}. Both $a$ and $b$ are incident to $v$. Let $a=vx$,
$b=vy$, $S=\{x,y\}$ and $\overline{S}=V(T_1)\backslash S$. Pick
$e\in[S,\overline{S}]_{T_1}$. Note that $e$ is not incident to
$v$. Similar to Case 1.1, $e$ satisfies (\ref{ineq1}).

Therefore, $E(T_1)$ is a resolving set of $L(G)$ with size
$|V(G)|-2$, and the upper bound is valid. $\qed$

The lower bound in Theorem \ref{upper} can be attained if $G$ is a
path. The fact that $\mu(L(K_{1,n}))=n-1$ implies that the upper
bound in Theorem \ref{upper} is tight. It seems to be difficult to
improve the bound for general graphs. However, for a tree $T$, we
can obtain the metric dimension of $L(T)$ in terms of some
parameters of $T$.

Let $T$ be a tree. A vertex of degree 1 in $T$ is called an
\emph{end-vertex}. A vertex of degree at least 3 in $T$ is called
a \emph{major vertex}. An end-vertex $u$ of $T$ is said to be a
\emph{terminal vertex of a major vertex} $v$ of $T$ if
$d_T(u,v)<d_T(u,w)$ for every other major vertex $w$ of $T$. A
major vertex $v$ of $T$ is an \emph{exterior major vertex} of $T$
if there exists a terminal vertex of $v$ in $T$. We denote the set
of all the exterior major vertices in $T$ by EX($T$); For
$v\in{\rm EX}(T)$, we denote the set of all the terminal vertices
of $v$ by TER($v$). Let $\sigma(T)=\sum_{v\in{\rm EX}(T)}|{\rm
TER}(v)|$ and ${\rm ex}(T)=|{\rm EX}(T)|$. Chartrand et al.
\cite{ChE} computed the metric dimension of a tree in terms of
$\sigma(T)$ and ${\rm ex}(T)$.

\begin{prop}\label{tree}{\em(\cite{ChE})}
If $T$ is a tree that is not a path, then $\mu(T)=\sigma(T)-{\rm
ex}(T)$.
\end{prop}

Finally, we shall compute the metric dimension of the line graph
of a tree. If $P$ is a path, then $\mu(L(P))=1$.

\begin{prop}\label{ltree}
If $T$ is a tree that is not a path, then
$\mu(L(T))=\sigma(T)-{\rm ex}(T)$.
\end{prop}
\proof Let $R$ be a resolving set of $L(T)$ with the minimum
cardinality. For a given vertex $v\in{\rm EX}(T)$, we claim that
\begin{equation}\label{lt1}
\sum_{u\in{\rm TER}(v)}|R\cap E(P_{u,v})|\geq |{\rm TER}(v)|-1,
\end{equation}
where $P_{u,v}$ is the unique path between $u$ and $v$ in $T$.
 To the contrary, suppose that there exist two
 different terminate vertices $u_1,u_2$ of $v$
 such that $R\cap E(P_{u_1,v})=R\cap E(P_{u_2,v})=\emptyset$.
 Let $e_1$ and $e_2$ be the edges incident to $v$ in $P_{u_1,v}$
 and $P_{u_2,v}$, respectively. For each $e\in R$, we have $d_{L(T)}(e_1,e)=d_{L(T)}(e_2,e)$,
 contradicting the fact that $R$ is a resolving set of $L(T)$.
 Hence our claim is valid.
 Since $|R|\geq\sum_{v\in{\rm EX}(T)}\sum_{u\in{\rm TER}(v)}|R\cap E(P_{u,v})|$,
 by (\ref{lt1}) we have
\begin{equation}\label{lt2}
\mu(L(T))=|R|\geq\sum_{v\in{\rm EX}(T)}(|{\rm
TER}(v)|-1)=\sigma(T)-{\rm ex}(T).
\end{equation}

Let $W$ be a set obtained from the end-vertex set of $T$ by
deleting one terminal vertex of each exterior major vertex of $T$.
 In \cite[Theorem 5]{ChE}, Chartrand et al. proved that $W$ is a
resolving set of $T$ with size $\sigma(T)-{\rm ex}(T)$. Let $W_L$
be the set of all the edges each of which is incident to one
vertex of $W$. Then $|W_L|=|W|$. We will show that $W_L$ is a
resolving set of $L(T)$.

For any two distinct edges $a$ and $b$ of $T$, there exists a
unique path
$$
(w_0,w_1,\ldots,w_{k-1},w_k)
$$
 such that $a=w_0w_1$
and $b=w_{k-1}w_k$. Since $w_0\neq w_k$, there exists a vertex
$w\in W$ such that $d_T(w_0,w)\neq d_T(w_k,w)$. Without loss of
generality, assume that $d_T(w_0,w)<d_T(w_k,w)$. Let $e$ be the
edge incident to $w$. Then $e\in W_L$.

{\em Case 1}. $w_1\in V(P_{w_0,w})$. Then
$$
d_{L(T)}(a,e)=d_T(w_0,w)-1<d_T(w_k,w)-1\leq d_{L(T)}(b,e).
$$

{\em Case 2}. $w_1\not\in V(P_{w_0,w})$.  Then
$(w_k,w_{k-1},\ldots,w_1,P_{w_0,w})$ is the unique path between
$w_k$ and $w$. It follows that
$$
d_{L(T)}(a,e)=d_T(w_0,w)<d_T(w_{k-1},w)= d_{L(T)}(b,e).
$$

Therefore, $W_L$ is a resolving set of $L(T)$, which implies that
$\mu(L(T))\leq\sigma(T)-{\rm ex}(T)$. By (\ref{lt2}), the desired
result follows. $\qed$

Combing Proposition \ref{tree} and Proposition \ref{ltree},
$\mu(T)=\mu(L(T))$ for a tree $T$. It seems to be interesting to
characterize a graph $G$ satisfying $\mu(G)=\mu(L(G))$.

\section*{Acknowledgement} This research is supported by NSF of China (10871027), NCET-08-0052,  and   the
Fundamental Research Funds for the Central Universities of China.

\end{CJK*}
\end{document}